\documentclass[11pt,twoside]{amsart}
\usepackage{amssymb,latexsym,amsthm,amsmath,mathtools,amsfonts,hyperref,fancyhdr,comment, tabularht, tabularx, longtable, mathrsfs, cleveref,booktabs, color, geometry,bm}
\usepackage[
singlelinecheck=false 
]{caption}

\usepackage{geometry}
\hypersetup{
	colorlinks=true,
	linkcolor=blue,
	filecolor=magenta,
	urlcolor=cyan,
	citecolor=red,
}

\long\def\symbolfootnote[#1]#2{\begingroup%
	\def\thefootnote{\fnsymbol{footnote}}\footnote[#1]{#2}\endgroup}

\newcommand{\Z}{\ensuremath{\mathcal{Z}}}

\def \R {{\mathbb{R}}}

\def \N {{\mathbb{N}}}
\def \Z {{\mathbb{Z}}}

\makeatletter
\def\imod#1{\allowbreak\mkern10mu({\operator@font mod}\,\,#1)}
\makeatother

\newtheorem{theorem}{Theorem}[section]
\newtheorem{lemma}[theorem]{Lemma}

\newtheorem*{theorem*}{Theorem}
\theoremstyle{definition}

\newtheorem{problem}[theorem]{Problem}

\numberwithin{equation}{section}
\newcommand{\ignore}[1]{}

\newcommand{\mynote}[1]{}

\usepackage[labelsep=period]{caption}
\begin{document}
	
	\title{Groups having 12 cyclic subgroups}
	
	\author{Khyati Sharma}
	\address{Shiv Nadar Institution of Eminence, NH-91, Dadri, Gautam Buddha Nagar}
	\email{khyatisharma0907@gmail.com}
	
	\author{A. satyanarayana Reddy}
	\address{Shiv Nadar Institution of Eminence, NH-91, Dadri, Gautam Buddha Nagar}
	\email{satya.a@snu.edu.in}
	
	\today 
	\begin{abstract}
		A finite group is said to be  {\em $n$-cyclic} if it contains $n$ cyclic subgroups. For a finite group $G$, the ratio of the number of cyclic subgroups to the number of subgroups is known as the {\em cyclicity degree} of the group $G$ and is denoted by $cdeg (G)$. In this paper, we classify all $12$-cyclic groups. We also prove that the set of cyclicity degrees for all the finite groups is dense in $[0,1]$, which gives a solution to the problem asked by T\u{a}rn\u{a}uceanu and T\'{o}th in \cite{CyclicityDegree} ``{\em For every $a\in [0, 1]$, does there exist a sequence $(G_n)_{n\in \N}$ of finite groups such that $\lim_{n\to\infty} cdeg(G_n)=a$} "?
	\end{abstract}
\maketitle
\section{Introduction}\label{sec1}
\par In the study of groups, subgroups play a fundamental role in understanding the structure of a group. By examining the number, types, and lattice diagram of subgroups, important structural features such as solvability, simplicity can be studied. Also, many enumeration problems in algebra and combinatorics often boil down to counting subgroups of certain related groups, either absolutely or up to some form of equivalence. A prime example is Galois theory: by enumerating subgroups (or certain types of
subgroups) of $S_n$, one gets precise information on the number of intermediate fields in a field extension of degree $n$. In this context, Betz and Nash~\cite{betz2022classifying} classified all abelian groups $G$ with $s(G)\leq 22$ and all non-abelian groups $G$ with $s(G)\leq 19$, where $s(G)$ denotes the number of subgroups of $G$. Additionally, Das and Mandal~\cite{das2024solvability} established certain criteria for the solvability and nilpotency of a finite group based on its number of subgroups. Recently, Dougal and Tracey~\cite{SymmetricGroupSubgroup} obtained a result concerning the number of subgroups in symmetric groups.\\
To study the subgroups of a finite group, it is natural to consider their cyclic subgroup structure, as every element of the group generates a cyclic subgroup. This counting is also helpful in group representation theory to find the Wedderburn decomposition of group algebras for an abelian group. We remark on various results related to the counting of cyclic subgroups. Miller was the first group theorist to give a result in this direction. In \cite{miller1929number}, Miller proved that: if $G$ is a finite group,  then
$$|G|=\sum\limits_{m||G|}c(m)\varphi (m),$$ where $c(m)$ denotes the number of cyclic subgroups of order $m$ in $G$. In another paper \cite{miller1905determination}, Miller proved that the number of cyclic subgroups of order $p^r, r>1$ in a finite group is always of the form $kp$ whenever the Sylow subgroups of the group are non-cyclic. Richard~\cite{richards1984remark} in $1984$ gave a remarkable lower bound on the number of cyclic subgroups. He showed that for a finite group $G$ of order $n$, the number of cyclic subgroups $c(G)$ satisfies the inequality $c(G)\geq \tau (n)$, where $\tau (n)$ denotes the number of divisors of $n$. Moreover, equality holds if and only if $G$ is a cyclic group of order $n$. Meng and Lu \cite{LowerBoundNumberofCyclicSubgroups} further refined this bound for finite non-cyclic nilpotent groups. T\'{o}th~\cite{c(G)for_abelian_groups} gave a formula to count the number of cyclic subgroups of finite abelian groups. In 2017, Jafari and Madadi~\cite{jafari2017number} also contributed to the counting of cyclic subgroups in finite groups.
\par Let $G$ be a finite group, and let $s(G)$ and $c(G)$ denote the number of subgroups and the number of cyclic subgroups of $G$, respectively. A group $G$ is {\em $n$-cyclic} if $c(G)=n$. It is easy to see that $c(G)\leq |G|$ and equality holds if and only if $G$ is an elementary abelian $2$-group. 
T\u{a}rn\u{a}uceanu~\cite{tuarnuauceanu2015finite} classified all finite groups $G$ having $|G|-1$ number of cyclic subgroups. In the same paper, the author posed a problem classifying the finite groups $G$ satisfying $c(G)=|G|-r$, where $2\leq r\leq (|G|-1)$. Motivated by his work, Belshoff et al.~\cite{belshoff2019finite} classified all finite groups having $|G|-r$ cyclic subgroups for $r=3,4$ and $5$. The authors further classified all finite groups having $|G|-r$ cyclic subgroups~ \cite{belshoff2018finite}, where $r$ is in the range $1-32$. Meanwhile, Zhou~\cite{zhou2016finite}, Kalra~\cite{kalra2019finite}, Ashrafi and Haghi~\cite{ashrafi2019n} classified all the $n$-cyclic groups for $3\leq n\leq 10$. In our earlier paper~\cite{11cyclic}, we classified all $11$-cyclic groups. 
\par Recently, many authors defined functions on the set of finite groups by using $c(G)$ and $s(G)$ and studied these functions and their relationship with the structure of group $G$. 
\par A new approach to study group-theoretic problems is to study them statistically and get probabilistic results, which help to understand
the asymptotic behavior. In this direction, the probability of a random subgroup of $G$ to be cyclic is introduced by T\u{a}rn\u{a}uceanu and T\'{o}th~\cite{CyclicityDegree}, known as the {\em cyclicity degree} of a finite group $G$ denoted by $cdeg (G)$ defined as $cdeg (G)=\frac{c(G)}{s(G)}$. The authors explicitly determine the cyclicity degree of some known classes of groups and study their various properties. In the same paper, the authors posed the following problem:
\begin{problem}\label{problem}
 For every $a\in [0, 1]$, does there exist a sequence $(G_n)_{n\in \N}$ of finite groups such that $\lim_{n\to\infty} cdeg(G_n)=a$?
 \end{problem}
\par In this article, we classify all $12$-cyclic groups. Additionally, we give a solution to the Problem~\ref{problem} posed by T\u{a}rn\u{a}uceanu and T\'{o}th~\cite{CyclicityDegree}, which is equivalent to prove that the image set of the function $cdeg$ is dense in $[0,1]$.\\
The organization of the article is as follows. In section~\ref{sec2}, we set up notations, recall the results
of \cite{11cyclic}, and formulate lemmas on which our proofs are based. The classification of $12$-cyclic groups and a solution of Problem~\ref{problem} are given in section~\ref{sec3} and section~\ref{sec4} respectively.
\section{Notations and Preliminaries}\label{sec2}
The notations $\Z_n$, $D_{2n}$, $Q_{2^n}$, $M(p^a)$, $F_5$, $SL(2,3)$ and $Dic_{n}$ denote the cyclic group of order $n$, dihedral group of order $2n$, generalized quaternion group of order $2^n$, modular group of order $p^a$, Frobenius group of order $20$, special linear group of order $24$ and dicyclic group of order $4n$ respectively.  The group $\Z_8.\Z_4=${\em SmallGroup$(32,15)$} has the following presentation $\langle a,b : a^8=1, b^4=a^4, bab^{-1}=a^{-1} \rangle$. Throughout this paper, $p, q$ and $r$ are distinct prime numbers, and the Euler totient function is denoted by $\varphi$. Here, $\N$ denotes the set of natural numbers and does not contain zero. The number of Sylow $p$-subgroups of $G$ is denoted by $n_p(G)$. A group $G$ is {\em CLT} if it has a subgroup corresponding to every divisor of $|G|$. A group $G$ is {\em Dedekind} if every subgroup of $G$ is normal. A non-cyclic group $G$ is said to be {\em minimal non-cyclic} if all its proper subgroups are cyclic.\\
Let $d(n), \omega(n)$ denote the number of positive divisors and the number of distinct prime divisors of $n$, respectively. Richard's Theorem~\cite{richards1984remark} provides a lower bound on $c(G),$ in particular $c(G)\geq d(|G|)$ and the equality holds if and only if $G$ is cyclic. Suppose $G$ is a $12$-cyclic group of order $n$. An immediate consequence of Richard's Theorem is that $\omega(n)\le 3$ and $n$ is one of the form $p^k, pq, p^2q, pqr, p^3q, p^4q, p^5q, p^2q^2, p^3q^2$ or $p^2qr$, where $k\leq {11}$. In this article we denote the quantity given in equation~\ref{equ2}, by $T(G)$ that is 
\begin{equation}\label{equ2}
    T(G)=|G|-\sum\limits_{m||G|}c(m)\varphi (m).
\end{equation}
In some cases, to find all $12$-cyclic groups, we take different possibilities of $c(m)$ for all the divisors $m$ of $|G|$ and then try to find all those groups $G$, which satisfies the equations $\sum\limits_{m||G|}c(m)=12$ and $T(G)=0$. In this paper, all the calculations are done by using GAP~\cite{GAP4} and the {\em SmallGroup$(n, i)$} denotes the $i^{th}$ group of order $n$ in the {\em Small Group Library} of GAP.
\begin{lemma}\label{lemma1}(Lemma 2.2,~\cite{11cyclic})  Let $G$ be a non-trivial $n$-cyclic group of order $p^aq^b$, where $p$ and $q$ are distinct prime numbers, and let a Sylow $p$-subgroup of $G$ denoted as $P$ is Dedekind. Then either $P$ is normal in $G$ or $q< {n-1}$.
    \end{lemma}
    \begin{lemma}\label{lemma2}(Proposition 2.4,~\cite{11cyclic}) Let $G$ be a non-cyclic, $n$-cyclic group, and let $M$ be a maximal subgroup of $G$. Then $c(M)<{n-1}$. 
        \end{lemma}
        \begin{lemma}\label{lemma3}(Frobenius,~\cite{frobenius1895verallgemeinerung})
    If $p$ is a prime, $p$ divides $|G|$, and $a_p$ denotes the number of subgroups of order $p$ in $G$, then $a_p\equiv 1\pmod p$.
\end{lemma}
\begin{lemma}\label{lemma4}(Lemma 4.1,~\cite{beta(G)}) Let $(x_n)_{n\geq 1}$ be a sequence of positive real numbers such that $\lim\limits_{n\to\infty}x_n=0$ and the series $\sum\limits_{i=1}^{\infty}x_i$ is divergent. Then, the set containing the sums of all finite subsequences of $(x_n)_{n\geq 1}$ is dense in $[0, \infty)$.
\end{lemma}
    
\section{Classification 12-cyclic groups}\label{sec3}
\begin{theorem}\label{main theorem1}
    Let $G$ be a finite group. Then $c(G)=12$ if and only if $G$ is isomorphic to an element of the set $S= \{\Z_{p^{11}}, \Z_{p^5q}, \Z_{p^3q^2}, \Z_{p^2qr}, {\Z_2}^2\times \Z_4, {\Z_2}^2\rtimes \Z_4, D_{16}, D_8\rtimes \Z_2, D_{18}, F_{5}, Dic_6, \Z_8.\Z_4, {\Z_{2}}^2\rtimes \Z_9, M(64), \Z_2\times \Z_{32}, \Z_5\times \Z_{25},  \Z_{25}\rtimes \Z_5, \Z_2\times \Z_{2s^2}, \Z_{4s}\times \Z_2 \}$, where $p, q, r$ are prime numbers and $s$ is an odd prime.
\end{theorem}
If $G$ is a $12$-cyclic group, then by Richard's Theorem~\cite{richards1984remark}, $|G|\in \{p^k, pq, p^2q, pqr, p^3q, p^4q, p^5q, \\ p^2q^2, p^3q^2, p^2qr\}$, where $k\leq {11}$. In addition, there is a unique $12$-cyclic group of orders $p^{11}, p^5q, p^3q^2$ and $p^2qr$, that is always cyclic. Moreover, for other possible orders of the $12$-cyclic groups, the proof is divided into the following four Lemmas.
\begin{lemma}\label{lem:proof1}
 If $G$ is a $12$-cyclic $p$-group, then
    $$G\cong \begin{cases}
      \Z_{p^{11}},\; {\Z_{2}}^2\times \Z_{4}\;, \Z_{2}\times \Z_{32}\; \mbox{or}\; \Z_{5}\times \Z_{25} & \mbox{if $G$ is abelian},\\
     {\Z_{2}}^2\rtimes \Z_{4}\;, D_8\rtimes \Z_2\;, D_{16}\;, \Z_8.\Z_4\;, M(64)\; \mbox{or}\; \Z_{25}\rtimes \Z_5 & \mbox{otherwise.}
    \end{cases}$$
    \end{lemma}
    \begin{proof}
    Let $M$ be a maximal subgroup of $G$. Then we have the following two cases:
\begin{enumerate}
 \item Let $M$ be a cyclic subgroup of $G$. We first assume that $G$ is abelian. Since $|M|=p^{a-1},$ then either $G$ is isomorphic to $\Z_{p^a}$ or $\Z_{p}\times \Z_{p^{a-1}}.$ If $G\cong \Z_{p^a}$, then it is easy to see that $a=11$. If $G\cong \Z_{p}\times \Z_{p^{a-1}}$, then by Theorem $1.1$ of \cite{ashrafi2019n}, $c(G)=(a-1)p+2,$ consequently, $p=5, a=3$ or $p=2, a=6$ and $|G|=125$ or $64.$ By using a simple GAP program~\cite{GAP4}, $G\cong \Z_{5}\times \Z_{25}=${\em SmallGroup$(125,2)$} or $\Z_{2}\times \Z_{32}=${\em SmallGroup$(64,50)$}.
 
 Now let $G$ be a non-abelian group, then depending on $p$, we have the following two situations. If $p$ is odd, then $G\cong M(p^a)$ and by Theorem $1.1$ \cite{ashrafi2019n}, $c(G)=(a-1)p+2.$ This shows that $p=5, a=3.$ By \cite{GAP4}, $G\cong \Z_{25}\rtimes \Z_5=${\em SmallGroup$(125,4)$}. If $p=2$, then by the classification theorem of finite non-abelian $2$-groups containing cyclic maximal subgroup, $G\cong D_{2^a}, Q_{2^a}, M(2^a)$ or $S_{2^a}.$ Using~\cite[Theorem $1.1$]{ashrafi2019n}, we have
        $c(D_{2^a})=a+2^{a-1},$ $c(Q_{2^a})=a+2^{a-2},$ $c(M(2^a))=2a$ and $c(S_{2^a})=a+3.2^{a-3}.$ Now by simple calculation, one can check that $G\cong M(64)$ or $D_{16}$.
\item Let $M$ be a non-cyclic subgroup of $G.$ By Lemma~\ref{lemma2}, $c(M)\leq 10$. Also by Theorem $1.2$ and $2.4$ \cite{ashrafi2019n}, $M$ is isomorphic to a member of the set $S=\{\Z_2\times \Z_2, \Z_3\times \Z_3, \Z_2\times \Z_4, \Z_5\times \Z_5, \Z_3\times \Z_9, \Z_2\times \Z_8, \Z_2\times \Z_2 \times \Z_2, \Z_7\times \Z_7, \Z_2\times \Z_{16}, \Z_4\times \Z_4,  \Z_2\times Q_{8}, SD_{16}, \Z_{16}\rtimes \Z_2, \Z_4\rtimes \Z_4, Z_9\rtimes \Z_3, \Z_8\rtimes \Z_2, Q_{16}, D_{8}, Q_8 \}$. Therefore we get the following possibilities, if $p=2$, then $a=3,4,5$ or $6,$ if $p=3$ then $a=3,4,$ or $5$ and if $p=5,7$ then $a=3.$ Moreover $|G|\in \{8, 16, 27, 32, 64, 81, 125, 243, 343\}.$ By using a simple GAP~\cite{GAP4} program one can see that $G\cong {\Z_2}^2\times \Z_4=${\em SmallGroup$(16,10)$}, ${\Z_2}^2\rtimes \Z_4=${\em SmallGroup$(16,3)$}, $D_8\rtimes \Z_2=${\em SmallGroup$(16,13)$}  and $\Z_8.\Z_4=$ {\em SmallGroup }{\em$(32,15)$}.
\end{enumerate}
\end{proof}
Hence by using the above Lemma and Richard's theorem~\cite{richards1984remark}, we conclude that a group $G$ is cyclic and $12$-cyclic if and only if $G\cong \Z_{p^{11}}, \Z_{p^5q}, \Z_{p^3q^2}, \Z_{p^2qr}$. Again, by using the above Lemma and \cite[Theorem 1]{c(G)for_abelian_groups}, it is easy to check that a group $G$ is non-cyclic abelian and $12$-cyclic if and only if $G\cong {\Z_2}^2\times \Z_4, \Z_2\times \Z_{32}, \Z_5\times \Z_{25}, \Z_2\times \Z_{2q^2}, \Z_{4q}\times \Z_2$, where $q$ is an odd prime number.\\
Therefore, from now onwards, all the groups are supposed to be non-abelian. The rest of the proof is organized as follows: In Lemma~\ref{lem:proof2}, we prove that if $G$ is a $12$-cyclic group of order $p^2q$ or $p^2q^2$, then $G$ is isomorphic to $F_5, D_{18}$ or ${\Z_2}^2\rtimes \Z_9$. In Lemma~\ref{lem:proof3}, we prove that if $G$ is a $12$-cyclic group of order $p^3q$, then $G$ is isomorphic to $Dic_6$. Further, in Lemma~\ref{lem:proof4}, we prove that there is no $12$-cyclic group of order $pq, pqr$ and $p^4q$.
\begin{lemma}\label{lem:proof2}
    Let $G$ be a $12$-cyclic group.
    \begin{enumerate}
        \item If $|G|=p^2q$, then $G$ is isomorphic to $F_5$ or $D_{18}$.
        \item If $|G|=p^2q^2$, then $G\cong {\Z_2}^2\rtimes \Z_9=$ { SmallGroup}$(36,3)$.
    \end{enumerate}
\end{lemma}
\begin{proof}
The proof contains the following cases.
    \begin{enumerate}
\item If $|G|=p^2q$, then there are the following two situations.
\begin{enumerate}
        \item If $p<q,$ then according to \cite[Proposition~$3.3$]{kalra2019finite}, $c(G)\in \{6, 2p+4, pq+4, q+4, 2q+2\}.$ Hence $G$ is $12$-cyclic only if  $p=2, 3$ and $q=5.$ Moreover, by \cite{GAP4} one can notice that $G\cong F_5=${\em SmallGroup$(20,3)$}.
        \item If $p>q$, then by~\cite[Proposition~$3.2$]{kalra2019finite},  $c(G)\in \{6, 2p+4, p^2+3, p^2+p+2, 2p+3, 3p+2\}.$ Thus $c(G)=12$ is possible only if $p=3$ and $q=2.$ Again by~\cite{GAP4}, we have $G\cong D_{18}=${\em SmallGroup$(18,1)$}. 
    \end{enumerate}
    \item If $|G|=p^2q^2$, then Sylow $p$ and $q$ subgroups of $G$ are Dedekind. Now, we prove the result by examining the following cases.
  
    \begin{enumerate}
    \item {\em \textbf{$\mathbf{G}$ has a unique subgroup of orders $\mathbf{p}$ and $\mathbf{q}$.}} In this case Sylow $p$ and $q$ subgroups of $G$ are cyclic. If $n_p(G)=n_q(G)=1$, then $G$ is cyclic, which is a contradiction. Also, by Sylow theorems, $n_p(G), n_q(G)>1$ is not possible. Let us assume that $n_p(G)>1$ and $n_q(G)=1$. Then by Sylow theorem, $n_p(G)\geq {1+p}$ and $n_p(G)\in \{q,q^2\}$. Moreover, $G$ has a unique cyclic subgroup of orders $1, p, q, pq$ and $q^2$.\\
            If $n_p(G)=q$, then $p<q$ and $c(G)\geq {q+5}$. Hence $p\in \{2,3,5\}$ and $q\in \{2,3,5,7\}$. Also $|G|\in \{36,100,196,225,441,1225\}$. If $n_p(G)=q^2$, then either $|G|=36$ or $c(G)>12$. A simple calculation with GAP~\cite{GAP4} shows that no such group of these orders is $12$-cyclic. 
            \item {\em \textbf{$\mathbf{G}$ has a unique subgroup of order $\mathbf{q}$ and at least $\mathbf{p+1}$ subgroups of order $\mathbf{p}$}.} This implies that Sylow $q$-subgroup of $G$ is cyclic. Therefore $c(G)\geq p+4,$ which implies that $p\in \{2,3,5,7\}.$ If Sylow {p}-subgroup of $G$ is not normal, then by Lemma~\ref{lemma1}, $q<11$ and $|G|\in \{36, 100, 196, 225, 441, 1225 \}$. Again if Sylow $p$-subgroup of $G$ is normal, then Sylow $q$-subgroup of $G$ is not normal (being $G$ non-abelian) that is $n_q(G)\geq {1+q}$. Since $G$ has at least $3$ subgroups of order $p$, then $c(G)\geq {6+q}$ and $q\leq 6$. Hence $|G|\in \{36, 100, 196, 225, 441, 1225\}$. By using a simple GAP~\cite{GAP4} program, one can check that $G\cong {\Z_2}^2\rtimes \Z_9=${\em SmallGroup$(36,3)$}.
            \item {\em \textbf{Subgroups of orders $\mathbf{p}$ and $\mathbf{q}$ are not unique.}} By Lemma~\ref{lemma3}, $G$ has at least $p+1$ and $q+1$ subgroups of order $p$ and $q$, respectively. Therefore, $c(G)\geq p+q+3$. Accordingly, $p+q\in \{2,3,4,5,6,7,8,9\}$ and $|G|\in\{36,100,196,225\}.$ By GAP~\cite{GAP4}, after analyzing 
 all the groups of order $36,100,196$ and $225$, it can be seen that none of them is $12$-cyclic.
\end{enumerate}
    \end{enumerate}
\end{proof}
\begin{lemma}\label{lem:proof3}
    If $G$ is a $12$-cyclic group of order $p^3q$, then $G\cong Dic_6$.
\end{lemma}
\begin{proof}
     Let $|G|=p^3q$. Then by using the proof of~\cite[Theorem 1.1]{11cyclic}, we can check that $G$ is a CLT group, so it has a subgroup of order $p^2q,$ let us call it $M$. By Lemma~\ref{lemma2}, we have $c(M)\leq 10$. Also by Theorem $1.2$ and $2.4$ \cite{ashrafi2019n}, $M\in \{\Z_{p^2q}, \Z_3\rtimes \Z_{4}, \Z_2\times \Z_{2p}, A_4, \Z_5\rtimes \Z_4,D_{12},\Z_3\times \Z_{3q} \}$. If $M\in \{A_4, \Z_3\rtimes \Z_{4},\Z_5\rtimes \Z_4, D_{12}\}$, then $|G|\in \{24,40 \}.$ It is easy to verify using GAP~\cite{GAP4}, that no such group of order $24$ and $40$ is $12$-cyclic. Now, we discuss the remaining possibilities of $M$ separately.
 \begin{enumerate}
     \item If $M\cong \Z_{p^2q}$, then $G$ has at least six cyclic subgroups of orders $1, p, q, p^2, pq$ and $p^2q$. Now, there are the following sub-cases:
 \begin{enumerate}
 \item {\em \textbf { $\mathbf{G}$ has a unique subgroup of order $\mathbf{p}$.}} In this case, Sylow $p$-subgroup of $G$ is either cyclic or generalized quaternion. First, suppose that $G$ has a unique cyclic Sylow $p$-subgroup, then $n_q(G)\geq 1+q$ and $n_q(G)\geq p$. This shows that $c(G)\geq 7+q$ and $c(G)\geq 6+p$. Thus $|G|\in \{24, 40, 54, 135, 250, 375\}$. It is easy to verify using GAP~\cite{GAP4} that no group of the above orders in which Sylow $p$-subgroup is cyclic and normal is $12$-cyclic. Therefore, if Sylow $p$-subgroup of $G$ is cyclic, then it is not normal, and $n_p(G)=q$. Moreover, $c(G)\geq 6+q$ and $q\in \{2, 3, 5\}$. Again by using Sylow theorem $n_p(G)\geq {1+p}$, which implies that $c(G)\geq 7+p$, and $p=2,3$ or $5$. Thus $|G|\in \{24, 40, 54, 135, 250, 375\}$. Now, with the help of a simple GAP program \cite{GAP4}, after checking the number of cyclic subgroups of all the groups of the orders obtained above, no such group is $12$-cyclic.\\
 If Sylow $p$-subgroup of $G$ is a generalized quaternion, then Sylow $p$ and $q$ subgroups of $G$ are Dedekind. Let us first assume that Sylow $q$-subgroup of $G$ is normal. If Sylow $p$-subgroup of $G$ is also normal, then $G\cong Q_8\times \Z_q$, and it is easy to see that $c(G)=10$. Hence, this case is not possible. Therefore Sylow $p$-subgroup of $G$ is not normal. Then by Lemma~\ref{lemma1}, $q<11$, which implies that $|G|\in \{24, 40, 56 \}$. After using a simple GAP~\cite{GAP4} program, no such group of these orders is $12$-cyclic. Now we are left with the case when Sylow $q$-subgroup of $G$ is not normal, then by Sylow theorem $n_q(G)=1+kq$, where $k\in \N$ and $n_q(G)\geq 4$. As a consequence, either $|G|=24$ or $c(G)>12$. By using a simple GAP~\cite{GAP4} program, one can check that $G\cong Dic_6=${\em SmallGroup$(24,4)$}.
  \item {\em \textbf{$\mathbf{G}$ has at least $\mathbf{p+1}$ subgroups of order $\mathbf{p}$.}} Since $c(M)=6$, then $c(G)\geq {6+p}$ and $p\leq 6$, now we discuss these possibilities separately.\\
     If $p=5$, then $G$ has exactly $6$ subgroups of order $5$. Now, there are the following two possibilities. If Sylow $5$-subgroup of $G$ is cyclic, then by Sylow theorem $c(G)>12$. If Sylow $5$-subgroup is not cyclic and contains an element of order $25$, then by using GAP~\cite{GAP4}, one can check that it always contains more than two cyclic subgroups of order $25$, which implies that $c(G)>12$. Hence, this case is not possible.\\
     If $p=3$, then by Lemma~\ref{lemma3}, $G$ has either $4$ or $7$ subgroups of order $3$. If $G$ has $7$ subgroups of order $3$, then by equation~\ref{equ2}, $T(G)=0$ has no solution. Thus from now onwards, assume that $G$ has $4$ subgroups of order $3$. By \cite[Table~1]{kalra2019finite}, Sylow $3$-subgroup of $G$ lies in $\{\Z_{27}, \Z_{9}\times \Z_3, \langle a,b|a^{9}=b^3=1, ba=a^{4}b\rangle\}$ as $G$ has an element of order $9$. If a Sylow $3$-subgroup $(P)$ of $G$ is isomorphic to $\Z_9\times \Z_3$, then it is easy to see that $c(P)=8$. Additionally, $G$ has at least one cyclic subgroup of orders $q, 3q$ and $9q$. Also, one can see that either Sylow $3$ or Sylow $q$ subgroup of $G$ is not normal, which implies that $c(G)>12$. If $P$ is isomorphic to $\langle a,b|a^{9}=b^3=1, ba=a^{4}b\rangle$, then by GAP~\cite{GAP4}, $c(P)=8$. This case is not possible, same as the previous case. Finally if $P$ is isomorphic to $\Z_{27}$, then $n_3(G)\geq 4$ and $c(G)>12$. Therefore, this case is not possible.\\
     If $p=2$, then by Lemma~\ref{lemma3}, $G$ has either $3, 5$ or $7$ subgroups of order $2$. Let us first assume that Sylow $q$-subgroup of $G$ is not normal, then $n_q(G)\geq {1+q}$ and $q=3$ by using the fact that $c(G)=12$. This shows that $|G|=24$, and we can verify by using GAP~\cite{GAP4} that no such group is $12$-cyclic. Therefore, from now onwards, Sylow $q$-subgroup of $G$ is normal. By \cite[Table~1]{kalra2019finite} Sylow $2$-subgroups of $G$ lie in the set $\{\Z_8, \Z_4\times \Z_2, Q_8, D_8\}$ as $G$ has a subgroup of order $4$. Now, there are the following two situations. At first, assume that Sylow $2$-subgroup of $G$ is normal, then $G\cong  Q_8\times \Z_q$ or $D_8\times \Z_q$, with some easy calculation, it can be seen that none of these groups is $12$-cyclic. Finally, assume that Sylow $2$-subgroup of $G$ is not normal. This further gives the following two possibilities. If Sylow $2$-subgroup of $G$ is Dedekind, then by Lemma~\ref{lemma1}, $q<11$. Thus $|G|\in \{24, 40, 56\}$, with the help of a simple program using GAP~\cite{GAP4}, we can check that no such group is $12$-cyclic. At last, we are left with the case when Sylow $2$-subgroup of $G$ is not Dedekind, then it is isomorphic to $D_8$. By Sylow theorem, $G$ has at least $3$ Sylow $2$ subgroups. Let $H\cong D_8$ and $K\cong D_8$ be the two Sylow $2$ subgroups of $G$. Then, by using the fact $|H\cap K|\leq 4$, one can see that $G$ has at least $7$ subgroups of order $2$. Now, by using equation~\ref{equ2}, there is no solution to the equation $T(G)=0$. Hence, no such group is $12$-cyclic.
 \end{enumerate}
\item  If $M\cong \Z_2\times \Z_{2q}$, then $c(M)=8$, $|G|=8q$ and $G$ has no cyclic subgroup of order $8$. Also, by Sylow theorem $n_2(G)\in \{1, q\}$ and $n_q(G)\in \{1,4,8\}$. If $n_q(G)=8$, then $c(G)>12$. Also, if $n_q(G)=4$, then $q=3$ and $|G|=24$. By using GAP~\cite{GAP4} no such group of order $24$ is $12$-cyclic.
Consequently, $G$ has a normal Sylow $q$-subgroup. Moreover, Sylow $2$-subgroup is neither isomorphic to $\Z_8$ nor $Q_8$ as $G$ has a subgroup isomorphic to $\Z_2\times \Z_2$. Further, if Sylow $2$-subgroup of $G$ is normal, then $G\cong D_8\times \Z_q$, and it is easy to check that $c(G)=14$. Hence, from now onwards, assume that Sylow $2$-subgroup is not normal. This gives two possibilities. First is Sylow $2$-subgroup is Dedekind, then by Lemma~\ref{lemma1}, $q<11$ and $|G|\in \{24, 40, 56 \}$. With the help of some simple program using GAP~\cite{GAP4}, we can check that no group of these orders satisfying the above conditions is $12$-cyclic. Finally, if Sylow $2$-subgroup is not Dedekind, then it is isomorphic to $D_8$. Also, $G$ has at least $3$ cyclic subgroups of order $2q$ and at least $8$ cyclic subgroups of order $1,2,4$ and $8$, then by equation~\ref{equ2}, $T(G)=0$ has no solution, which is a contradiction.
\item  If $M\cong \Z_3\times \Z_{3q}$, then $c(M)=10$, $|G|=27q$. By Sylow theorem, $n_3(G)\in \{1,q\}$ and $n_q(G)\in \{1,3,9,27\}$. If $n_q(G)\in \{3, 9, 27\}$, then we can easily check that either $c(G)>12$ or $G$ is not $12$-cyclic by using equation~\ref{equ2}. Therefore, $G$ has a normal Sylow $q$-subgroup. Also, by Theorem~$1.2$~\cite{ashrafi2019n} Sylow $3$-subgroup $(P)$ of $G$ is either isomorphic to $\Z_9\times \Z_3$ or $\Z_9\rtimes \Z_3$ and $c(P)=8$. This implies that $c(G)>12$, as $c(M)=10$. Thus, no such group is $12$-cyclic.
\end{enumerate}
\end{proof}
\begin{lemma}\label{lem:proof4}
    If $|G|=pq, pqr$ or $p^4q$, then $G$ is not $12$-cyclic.
\end{lemma}
\begin{proof}
We will discuss the proof casewise.
    \begin{enumerate}
        \item If $|G|=pq$, where $p<q$, then by \cite[Lemma~$3.1$]{kalra2019finite} $c(G)=4$ or $q+2$. Hence, no group of order $pq$ is $12$-cyclic.
        \item Let $|G|=pqr$. Since every group of square-free orders is solvable, then $G$ has Hall subgroups of orders $pq,pr$ and $qr$. These Hall subgroups are either cyclic or they are isomorphic to $S_3, D_{14}, D_{10}$ or $\Z_7\rtimes \Z_3$ by Theorem~$1.2$ and $2.4$ of \cite{ashrafi2019n}. If all the Hall subgroups of $G$ are cyclic, then $G$ is a minimal non-cyclic group, which is not possible by \cite[Proposition~2.8]{jafari2017number}. Let $M$ be a non-abelian, maximal Hall subgroup of $G$. Then we have the following sub-cases:
    \begin{enumerate}
     \item If $M\cong \Z_7\rtimes \Z_3$, then $|G|=21p$. One can observe that $G$ has a unique cyclic subgroup of orders $1$ and $7$, and at least $7$ cyclic subgroups of order $3$. By using Sylow theorems, Lemma~\ref{lemma3}, and $c(G)=12$, all the possibilities for the number of cyclic subgroups of $G$ along with the function $T(G)$ using equation~\ref{equ2} are recorded in Table~\ref{tab:pqr1}. Now, it is easy to check that $T(G)=0$ has no solution. Hence, no such group is $12$-cyclic.
     \begin{table}[ht!]
\centering
\begin{tabular}{|cccccccc|}
 \hline
 $c(1)$ & $c(7)$ & $c(3)$ & $c(p)$ & $c(7p)$ & $c(3p)$ & $c(21)$ & $T(G)$\\
 \hline
  $1$ & $1$ & $7$ & $3$ & $0$ & $0$ & $0$ & $p-1$\\
  $1$ & $1$ & $7$ & $1$ & $1$ & $1$ & $0$ & $p-1$\\
  $1$ & $1$ & $7$ & $1$ & $1$ & $0$ & $1$ & $7p-13$\\
  \hline
 \end{tabular}\\
 \caption{}
           \label{tab:pqr1}
       \end{table}
 \item If $M\cong D_{14}$, then $|G|=14p$ and $c(M)=9$ by Theorem~$2.4$~\cite{ashrafi2019n}. By Sylow theorem, there is a unique subgroup of orders $7$ and $p$. Also, $G$ has at least $7$ subgroups of order $2$ and a unique cyclic subgroup of order $7p$. By using the fact $c(G)=12$, either $G$ has a cyclic subgroup of order $2p$ or order $14$. In both the cases, by equation~\ref{equ2}, $T(G)=0$ has no solution. Therefore, no such group is $12$-cyclic.
 \item If $M\cong S_3$, then $|G|=2\cdot 3\cdot p$. By Sylow theorem $n_p(G)$ is either $1$ or $6$. If $n_r(G)=1$, then $M$ is not normal. In this case either $n_2(G)$ or $n_3(G)$ is $p$. Therefore $p\leq 7$ and $|G|=30$ or $42$. Also, if $n_p(G)=6$, then $|G|=30$. By GAP~\cite{GAP4} no such group of order $30$ and $42$ is $12$-cyclic.
 \item If $M\cong D_{10}$, then $|G|=2\cdot 5\cdot p$. By Sylow theorem $n_p(G)$ is either $1$ or $4$. If $n_p(G)=1$, then $M$ is not normal. This implies that either $n_2(G)$ or $n_5(G)$ is at least $p$. Therefore $p\leq 7$ and $G=30$ or $70$. Also, if $n_p(G)=4$, then $|G|=30$. By using GAP~\cite{GAP4} no such group of order $30$ and $70$ is $12$-cyclic.
 \end{enumerate}
 \item Let $|G|=p^4q$. First, suppose that $G$ is a CLT group. Then $G$ has a subgroup of order $p^3q$, say $M$ such that $c(M)\leq 10$ by Lemma~\ref{lemma2}. By Theorem~$1.2$ and $2.4$ of \cite{ashrafi2019n}, $M\in \{\Z_{p^3q}, \Z_3\rtimes \Z_8, Q_8\times \Z_q, \Z_5\rtimes \Z_8 \}$. If $M\cong \Z_3\rtimes \Z_8$ or $\Z_5\rtimes \Z_8$ then $|G|=48$ or $80$. By GAP~\cite{GAP4}, no group of these orders containing a maximal subgroup isomorphic to $\Z_3\rtimes \Z_8$ or $\Z_5\rtimes \Z_8$ is $12$-cyclic. If $M\cong Q_8\times \Z_q$, then $c(M)=10$, and Sylow $2$-subgroup of $G$ is non-abelian. By \cite[Table~2]{kalra2019finite}, we can check that $c(G)>12$. Finally, we are left with the case when $M\cong \Z_{p^3q}$ and $c(M)=8$. Now, we proceed further by discussing the following sub-cases.
 \begin{enumerate}
     \item {\em \textbf{ $\mathbf{G}$ has a unique subgroup of order $\mathbf{p}$.}} In this case, Sylow $p$-subgroup of $G$ is either cyclic or generalized quaternion. First, suppose that $G$ has a unique cyclic Sylow $p$-subgroup. Then by Sylow theorem $n_q(G)\geq {1+q}$ and $n_q(G)\geq p$, otherwise $G$ will be cyclic. Consequently, by using the fact $c(M)=8$ and $c(G)=12$, we get $p,q\in \{2, 3\}$ and $|G|=48$ or $162$. By using a simple GAP~\cite{GAP4} program, we can check that no such group of these orders is $12$-cyclic. Therefore, if Sylow $p$-subgroup of $G$ is cyclic, then it is not normal that is $n_p(G)\geq {1+p}$ and $n_p(G)=q$. This implies that $c(G)\geq {9+p}$, and $p=2$ or $3$. Also $c(G)\geq {q+8}$, and $q=2$ or $3$. Thus $|G|=48$ and $162$. Hence, this case is not possible by using GAP~\cite{GAP4}.
     \par \noindent If Sylow $p$-subgroup of $G$ is generalized quaternion $(Q_{16})$, then $|G|=16q$ and $c(Q_{16})=8$. Also, $G$ has a unique cyclic subgroup of orders $q, 2q, 4q$ and $8q$. Now, by using equation~\ref{equ2}, it is easy to see that $T(G)=0$ has no solution.
 \item {\em \textbf{ $\mathbf{G}$ has at least $\mathbf{p+1}$ subgroups of order $\mathbf{p}$.}} In this case, by using the fact $c(G)=12$ it is easy to see that either $p=2$ or $3$. If Sylow $q$-subgroup is not normal, then $n_q(G)={1+kq}$, where $k\in \N$. Therefore $c(G)\geq {10+q}$ as $c(M)=8$, which implies that $q=2$ and $|G|=162$. By using GAP~\cite{GAP4}, no such group of order $162$ is $12$-cyclic. Now we are left with the case when Sylow $q$-subgroup is normal. In this case, by using \cite[Table 2,3]{kalra2019finite} any Sylow $p$-subgroup $(P)$ lies in the set $\{\Z_{p^{4}}, \Z_8\times \Z_2, Q_{16}, M(16)\}$ as $c(P)<{12}$. If $P$ is isomorphic to $\Z_{p^{4}}$, then $n_p(G)\geq 3$, which implies that $c(G)>{12}$. If $P\cong Q_{16}$, then $G$ has at least five subgroups of order $4$. Moreover, $G$ has at least three subgroups of order $2$ and at least one cyclic subgroup of orders $8, q, 2q, 4q$ and $8q$. Therefore $c(G)> 12$. Also, if Sylow $p$-subgroup is isomorphic to $\Z_8\times \Z_2$ or $M(16)$, then by using equation~\ref{equ2}, $T(G)=0$ has no solution.\\ 
 Hence, there is no CLT group of order $p^4q$, which is $12$-cyclic.
 \end{enumerate}
 Let $G$ be a non-CLT group of order $p^4q$. Then by \cite[Table 2,3]{kalra2019finite} any Sylow $p$-subgroup $P$ lies in the set $\{\Z_{p^{4}}, \Z_8\times \Z_2, \Z_4\times \Z_4, \Z_{27}\times \Z_3, G_7, G_{10}, G_{11}, G_{13}, G_{14}, G(vi)\}$. Now, we discuss these cases separately.\\
  If $P$ is isomorphic to $\Z_{27}\times \Z_3$ or $G(vi)\cong \Z_{27}\rtimes \Z_3$ then by \cite[Table 2,3]{kalra2019finite}, $c(G)=11$. Also, $G$ has a unique subgroup of order $q$. This implies that $|G|=80+q$, which is a contradiction.\\
 If $P$ is isomorphic to $\Z_4\times \Z_4, G_7, G_{10}, G_{13}$, $\Z_8\times \Z_2, G_{14}\cong Q_{16}, G_{11}\cong M(16)$, then by \cite[Table 2]{kalra2019finite} we get $c(P)=8$ or $10$. In this case, if Sylow $q$-subgroup is not normal, then by Sylow theorems $n_q(G)\geq 4$, which shows that $c(G)>12$ or $|G|=48$. By using GAP~\cite{GAP4}, no such group of order $48$ is $12$-cyclic. Therefore, from now onwards $G$ has a normal Sylow $q$-subgroup. Now, one can easily check that $G$ is a CLT group. Hence, this case is not possible.\\
 If $P\cong \Z_{p^4}$, then $c(P)=5$. In this case, $G$ does not contain any normal Sylow subgroup. Otherwise, $G$ will be a CLT group. Then by using Sylow theorems and the fact $c(G)=12$, we get $p, q\in \{2,3\}$ and $|G|= 48$ or $162$. By using a simple GAP~\cite{GAP4} program, we can check that no such group of above orders is $12$-cyclic.
    \end{enumerate}
\end{proof}
\textbf{Proof of theorem}~\ref{main theorem1} - 
 The proof follows from Lemmas~\ref{lem:proof1},\ref{lem:proof2},\ref{lem:proof3} and \ref{lem:proof4}.
    \section{Solution to the problem~\ref{problem}}\label{sec4}
    We denote the class of all finite groups by $\mathscr{F}$. The cyclicity degree $(cdeg(G)=\frac{c(G)}{s(G)})$ measures the probability of a random subgroup of $G$ to be cyclic. It satisfies the relation $0<cdeg(G)\leq 1$, and equality holds if and only if $G$ is cyclic. Also $cdeg$ is a function from $\mathscr{F}$ to $[0,1]$ with image set $S=\{cdeg(G) | G\in \mathscr{F}\}\subseteq (0, 1]$.\\
    Now we give a solution to the Problem~\ref{problem} by restating it as follows.
    \begin{theorem}\label{main theorem2}
    The set $\{cdeg(G) | G\in \mathscr{F}\}$ is dense in $[0, 1]$, where $\mathscr{F}$ denotes the class of all finite groups.
\end{theorem}
 \begin{proof} 
Let $p_n$ denotes the $n^{th}$ prime number, where $n$ is a positive integer. Then we have
    $$cdeg(\Z_{p_n}\times \Z_{p_n})=\frac{p_n+2}{p_n+3}.$$
            Consider the sequence $(x_n)_{n\geq 1}\subset(0,\infty),$ where $x_n=\ln{(\frac{p_n+3}{p_n+2})}.$ We have
        \[\lim_{n\to\infty}\frac{x_n}{\frac{1}{p_n}}=\lim_{n\to\infty}\frac{\ln{(\frac{p_n+3}{p_n+2}})}{\frac{1}{p_n}}=1.\]
        Since the series $\sum\limits_{n=1}^{\infty}\frac{1}{p_n}$ is divergent, then by \cite[Theorem 1]{hoang2015limit}, the series $\sum\limits_{n=1}^{\infty}x_n$ is also divergent. Also, $\lim\limits_{n\to\infty}x_n=0$ so all the hypotheses of Lemma~\ref{lemma4} are satisfied.\\ 
        Therefore we have \[\overline{\bigg\{\sum_{i\in I}\ln{\bigg(\frac{p_i+3}{p_i+2}\bigg)}|i\subset \N, |I|<\infty, \text{where $p_i$ is the $i^{th}$ prime number}\bigg\}}=[0,\infty)\Longleftrightarrow\]
        \[\overline{\bigg\{\ln{\bigg(\prod_{i\in I}\frac{p_i+3}{p_i+2}\bigg)}|i\subset \N, |I|<\infty, \text{where $p_i$ is the $i^{th}$ prime number}}\bigg\}=[0,\infty).\]
        Further we denote the interval $(0,\infty)$ by $Y$. Consider the topological spaces $(\R, \tau_{\R})$ and $(Y,\tau_{Y})$, where $\tau_{\R}$ is the usual topology of $\R$ and $\tau_{Y}$ is the subspace topology on $Y$. Since the function
        \[exp:(\R, \tau_{\R})\to (\R, \tau_{\R}), \text{given by $exp(x)=e^x, \forall x\in \R$},\] is continuous and $[1,\infty)$ is closed set of $\R$, we have  \[\overline{\bigg\{\prod_{i\in I}\frac{p_i+3}{p_i+2}|i\subset \N, |I|<\infty, \text{where $p_i$ is the $i^{th}$ prime number}}\bigg\}=[1,\infty).\] Note that
        \[\overline{\bigg\{\prod_{i\in I}\frac{p_i+3}{p_i+2}|i\subset \N, |I|<\infty, \text{where $p_i$ is the $i^{th}$ prime number}}\bigg\}_{\tau_Y}=\]
        \[\overline{\bigg\{\prod_{i\in I}\frac{p_i+3}{p_i+2}|i\subset \N, |I|<\infty, \text{where $p_i$ is the $i^{th}$ prime number}}\bigg\}\cap Y=\]
        \[[1,\infty)\cap Y=\overline{[1,\infty)}\cap Y=\overline{[1,\infty)}_{\tau_Y}.\] Hence, if we consider the continuous function \[f:Y \to \R, \text{given by $f(y)=\frac{1}{y}, \forall y\in Y$},\] we deduce that \[\overline{\bigg\{\prod_{i\in I}\frac{p_i+2}{p_i+3}|i\subset \N, |I|<\infty, \text{where $p_i$ is the $i^{th}$ prime number}}\bigg\}=\overline{(0,1]}=[0,1].\] Also
        \[\overline{\bigg\{\prod_{i\in I}cdeg(\Z_{p_i}\times \Z_{p_i})|i\subset \N, |I|<\infty, \text{where $p_i$ is the $i^{th}$ prime number}}\bigg\}=[0,1].\]
        By using \cite[Proposition 2.2]{CyclicityDegree} above set can be rewritten as
        \[\overline{\bigg\{cdeg (\bigtimes_{i\in I} (\Z_{p_i}\times \Z_{p_i})),i\subset \N, |I|<\infty, \text{where $p_i$ is the $i^{th}$ prime number}}\bigg\}=[0,1].\] Also we have \[\bigg\{cdeg (\bigtimes_{i\in I} (\Z_{p_i}\times \Z_{p_i})),i\subset \N, |I|<\infty, \text{where $p_i$ is the $i^{th}$ prime number}\bigg\}\subseteq S \subseteq (0,1].\] Then, by taking closures of the above three sets, we deduce that the set $S$ is dense in $[0,1]$. This completes the proof. 
         \end{proof}  
          
   \noindent \textbf{Acknowledgement:} The first-named author is supported by the University Grant Commission (UGC), India, under the scheme UGC-SRF.\\

\noindent \textbf{Conflict of interest:} The authors declare that they have no known competing financial interests or personal
relationships that could have appeared to influence the work reported in this paper.\\
  \bibliographystyle{abbrv}
	\bibliography{refs.bib} 
\end{document}